\documentclass[reqno]{amsart}

\makeatletter
\@addtoreset{equation}{section}

\makeatother

\newcommand{\Rm}{\mathbb{R}}

\newcommand{\be}{\begin{equation}}
\newcommand{\ee}{\end{equation}}
\newcommand{\ba}{\begin{equation}\begin{aligned}}
\newcommand{\ea}{\end{aligned}\end{equation}}

\newcommand{\pp}{\partial}
\newcommand{\vv}[1]{\boldsymbol{\mathrm{#1}}}

\newcommand{\argmin}{\mathop{\mathrm{arg\,min}}}

\usepackage{amsmath}
\usepackage[foot]{amsaddr}
\usepackage{amsthm,amssymb,mathrsfs}
\usepackage{graphicx}
\usepackage{color}

\theoremstyle{remark}

\title[]{Diffuse optical tomography by simulated annealing via a spin Hamiltonian}

\author{Yu Jiang$^1$}
\author{Manabu Machida$^2$}
\author{Norikazu Todoroki$^3$}
\address{${^1}$School of Mathematics, Shanghai University of Finance and Economics, Shanghai 200433, P.R. China}
\address{${^2}$Institute for Medical Photonics Research, Hamamatsu University School of Medicine, Hamamatsu 431-3192, Japan}
\address{${^3}$Department of Physics, Chiba Institute of Technology, Chiba 275-0023, Japan}
\email{jiang.yu@mail.shufe.edu.cn}
\email{machida@hama-med.ac.jp}
\email{todoroki.norikazu@p.chibakoudai.jp}

\date{\today}

\begin{document}

\begin{abstract}
Diffuse optical tomography (DOT) is an imaging modality which uses near-infrared light. Although iterative numerical schemes are commonly used for its inverse problem, correct solutions are not obtained unless good initial guesses are chosen. We propose a numerical scheme of DOT which works even when good initial guesses of optical parameters are not available. We use simulated annealing (SA) which is a method of the Markov-chain Monte Carlo. To implement SA for DOT, a spin Hamiltonian is introduced in the cost function, and the Metropolis algorithm or single-component Metropolis-Hastings algorithm is used. By numerical experiments, it is shown that an initial random spin configuration is brought to a converged configuration by SA and targets in the medium are reconstructed. The proposed numerical method solves the inverse problem for DOT by finding the ground state of a spin Hamiltonian with SA.
\end{abstract}

\maketitle

\section{Introduction}

In near-infrared spectroscopy, tomographic images of optical properties are obtained by diffuse optical tomography (DOT) \cite{Gibson-etal05,Yamada-Okawa14}. To obtain reconstructed images, inverse problems of determining coefficients of the diffusion equation from boundary measurements are solved.

For these inverse problems in DOT, direct approaches such as the use of the Born approximation and iterative methods such as the conjugate gradient method and damped Gauss-Newton method are commonly used \cite{Arridge99,Arridge-Schotland09}. One of the direct approaches is the inversion of a linearized problem by singular value decomposition (SVD) with the $L^2$ regularization (for example, see \cite{Markel-Schotland04}). To solve the nonlinear inverse problem of DOT, iterative methods are frequently used \cite{Arridge99}. For example, breast cancer was detected with the nonlinear-conjugate gradient method \cite{Choe-etal05}. See, e.g., \cite{Zhao-Cooper17,Zhu-Poplack20} for recent advances in DOT.

Although different penalty terms can be employed in iterative methods, the calculation gets trapped by local minima unless a good initial guess is given (for example, see \cite{Schweiger-Arridge99}). It was demonstrated in \cite{Jiang-etal19} that bad initial guesses (even though they look reasonable) lead completely wrong reconstructed results even when the diffusion coefficient is known and the absorption coefficient is characterized by only one unknown parameter in the diffusion equation.

As numerical schemes which are free from the issue of being trapped by local minima, statistical approaches have been developed \cite{Arridge-etal06,Shimokawa-etal12}. Monte Carlo methods have the advantage that the computation is not trapped by local minima since jumps occur randomly in the landscape. In particular, the Metropolis-Hastings Monte Carlo algorithm was used to estimate parameters in coefficients of the diffusion equation or the radiative transport equation, which is approximated to the diffusion equation at large scales \cite{Bal-etal13,Langmore13,Jiang-etal19}. However, only several unknown parameters such as the position of a target, and the absorption and scattering coefficients of the target, were determined in those studies.

In this paper, values of the absorption coefficient will be determined at $1830$ different points in the medium. The simulated annealing (SA) (e.g., \cite{Landau-Binder}) makes it possible to handle a large number of unknown parameters. To our knowledge, this is the first attempt to use SA for DOT.

The rest of the paper is organized as follows. In Sec.~\ref{DE}, we introduce the diffusion equation, which near-infrared light in biological tissue obeys. In Sec.~\ref{rytov}, we formulate our inverse problem with the Rytov approximation. Spins are introduced in Sec.~\ref{discr}. Then the cost function for the inverse problem is rewritten as a spin Hamiltonian in Secs.~\ref{spinh1} and \ref{spinh}. The algorithm of SA for DOT is described in Sec.~\ref{SA}. In Sec.~\ref{numerics}, numerical tests of our SA are presented. Finally, concluding remarks are given in Sec.~\ref{concl}.

\section{Diffusion equation}
\label{DE}

We consider DOT in a domain $\Omega\subset\Rm^d$ ($d=2,3$). Let $\pp\Omega$ be the boundary of $\Omega$. Let $\vv{\nu}(\vv{r})$ be the outer unit normal vector at $\vv{r}\in\pp\Omega$. We assume that $\Omega$ is occupied by biological tissue and the outside of $\Omega$ is vacuum. Let $u$ be the photon density of near-infrared light.

If the reconstruction is done using time-resolved data, the time-dependent diffusion equation needs to be considered. Let $c$ be the speed of light in $\Omega$. We assume that the diffusion coefficient $D_0$ is a positive constant but $\mu_a$ varies in space. The energy density $u$ obeys
\be
\left\{\begin{aligned}
&
\left(\frac{1}{c}\frac{\pp}{\pp t}-D_0\Delta+\mu_a(\vv{r})\right)
u(\vv{r},t)=f(\vv{r},t),
\quad\vv{r}\in\Omega,\quad t>0,
\\
&
D_0\vv{\nu}\cdot\nabla u(\vv{r},t)+\frac{1}{\zeta}u(\vv{r},t)=0,
\quad\vv{r}\in\pp\Omega,\quad t>0,
\\
&
u(\vv{r},0)=0,\quad\vv{r}\in\Omega,
\end{aligned}\right.
\label{de2d}
\ee
where $\zeta>0$ is a constant (see Sec.~\ref{numerics}) and $f(\vv{r},t)$ is the source term. We note that $\mu_a$ is nonnegative.

In the case of the continuous-wave measurement, light propagation in $\Omega$ is governed by the following diffusion equation for $u(\vv{r})$ ($\vv{r}\in\Omega$).
\be
\left\{\begin{aligned}
-D_0\Delta u(\vv{r})+\mu_a(\vv{r})u(\vv{r})=f(\vv{r}),
\quad \vv{r}\in\Omega,
\\
D_0\vv{\nu}\cdot\nabla u(\vv{r})+\frac{1}{\zeta}u(\vv{r})=0,
\quad \vv{r}\in\pp\Omega.
\end{aligned}\right.
\label{de0}
\ee
Even for measurements in time domain, quite often the time-independent diffusion equation (\ref{de0}) is used for reconstruction with moments or Laplace transforms of the raw time-resolved data. We assume the incident beam $f(\vv{r})$ as
\be
f(\vv{r},t)=g_0\delta(\vv{r}-\vv{r}_s^{(p)})\delta(t)
\quad\mbox{or}\quad
f(\vv{r})=g_0\delta(\vv{r}-\vv{r}_s^{(p)}),
\ee
where $g_0>0$ is a constant and $\vv{r}_s^{(p)}$ is the position of the source of the $p$th source-detector pair ($p=1,2,\dots,M_{\rm SD}$). Furthermore, $\delta(\cdot)$ is the Dirac delta function. Light is detected on the boundary at $\vv{r}_d^{(p)}$.

We choose a region of interest $\Omega_{\rm ROI}$ in $\Omega$. Let $\bar{\mu}_a$ be a nonnegative constant. We suppose
\be
\mu_a(\vv{r})=\bar{\mu}_a
\ee
for $\vv{r}$ in $\Omega$ including the boundary $\pp\Omega$ and excluding $\Omega_{\rm ROI}$. Let us write $\mu_a(\vv{r})$ as
\be
\mu_a(\vv{r})=\bar{\mu}_a+\delta\mu_a(\vv{r}),\quad \vv{r}\in\Omega_{\rm ROI}.
\ee

\section{Rytov approximation}
\label{rytov}

We will use SA by regarding $\delta\mu_a(\vv{r})$ as spins. For this purpose, here we develop the perturbation theory and give the solution $u(\vv{r})$ in terms of an integral of $\delta\mu_a(\vv{r})$. Since the Rytov approximation for the time-dependent case is similarly derived, we consider the time-independent case below.

Let $u_0(\vv{r})$ be the solution to (\ref{de0}) in the case of $\delta\mu_a\equiv0$. Then $u_0(\vv{r})=g_0G(\vv{r},\vv{r}_s^{(p)})$, where $G(\vv{r},\vv{r}')$ is the Green's function, which satisfies (\ref{de0}) with the source term $f(x)=\delta(\vv{r}-\vv{r}')$. Let us subtract the equation for $u_0$ from the equation for $u$ \cite{Arridge99}:
\be
\left\{\begin{aligned}
-D_0\Delta \left(u(\vv{r})-u_0(\vv{r})\right)+\bar{\mu}_a\left(u(\vv{r})-u_0(\vv{r})\right)=-\delta\mu_a(\vv{r})u(\vv{r}),
\quad \vv{r}\in\Omega,
\\
D_0\vv{\nu}\cdot\nabla \left(u(\vv{r})-u_0(\vv{r})\right)+\frac{1}{\zeta}\left(u(\vv{r})-u_0(\vv{r})\right)=0,
\quad \vv{r}\in\pp\Omega.
\end{aligned}\right.
\ee
Since the above diffusion equation for $u-u_0$ has the source term $-\delta\mu_au$, the solution $u(\vv{r})$ satisfies the following identity.
\be
u(\vv{r})=u_0(\vv{r})
-\int_{\Omega}G(\vv{r},\vv{r}')\delta\mu_a(\vv{r}')u(\vv{r}')\,d\vv{r}'.
\ee
The $n$th Born approximation $u_n$ is given by $u_n=\sum_{k=0}^nv_k$, where
\be
v_{k+1}(\vv{r})=-\int_{\Omega}G(\vv{r},\vv{r}')\delta\mu_a(\vv{r}')v_k(\vv{r}')\,d\vv{r}',\quad k=0,1,\dots,
\ee
with $v_0(\vv{r})=u_0(\vv{r})$. When the perturbation $\delta\mu_a$ is small, the Born series converges and $u=\lim_{n\to\infty}u_n$. In particular for nonnegative $\delta\mu_a$, a concise proof of the convergence is known under the condition $0\le\delta\mu_a(\vv{r})\le\bar{\mu}_a$ for any $\vv{r}\in\Omega$ \cite{Markel-Schotland07}. Let us define the first and second Rytov approximations $u_R,u_{R2}$ as
\be
u_R=u_0e^{v_1/u_0},\quad
u_{R2}=u_0e^{v_1/u_0}\exp\left[\frac{v_2}{u_0}-\frac{1}{2}\left(\frac{v_1}{u_0}\right)^2\right].
\ee
When the first Born approximation is compared with the first Rytov approximation, the superiority of the latter has been discussed \cite{Arridge99,Keller69,Kirkinis08}.

The difference $u-u_R$ can be written as
\ba
u(\vv{r})-u_R(\vv{r})
&=
u_0(\vv{r})\left(1-e^{v_1(\vv{r})/u_0(\vv{r})}\right)
\\
&-
\int_{\Omega}G(\vv{r},\vv{r}')\delta\mu_a(\vv{r}')u_R(\vv{r}')\,d\vv{r}'
\\
&-
\int_{\Omega}G(\vv{r},\vv{r}')\delta\mu_a(\vv{r}')\left(u(\vv{r}')-u_R(\vv{r}')\right)\,d\vv{r}'.
\ea
When $\|\mu_a\|_{L^{\infty}(\Omega)}$ is small, we have
\be
1-e^{v_1(\vv{r})/u_0(\vv{r})}\approx
-\frac{v_1(\vv{r})}{u_0(\vv{r})}=
\frac{1}{u_0(\vv{r})}\int_{\Omega}G(\vv{r},\vv{r}')\delta\mu_a(\vv{r}')u_0(\vv{r}')\,d\vv{r}',
\ee
and
\ba
\int_{\Omega}G(\vv{r},\vv{r}')\delta\mu_a(\vv{r}')u_R(\vv{r}')\,d\vv{r}'
&=
\int_{\Omega}G(\vv{r},\vv{r}')\delta\mu_a(\vv{r}')u_0(\vv{r}')e^{v_1(\vv{r}')/u_0(\vv{r}')}\,d\vv{r}'
\\
&\approx
\int_{\Omega}G(\vv{r},\vv{r}')\delta\mu_a(\vv{r}')u_0(\vv{r}')\,d\vv{r}'.
\ea
Therefore for sufficiently small $\|\delta\mu_a\|_{L^{\infty}(\Omega)}$, we obtain
\be
\|u-u_R\|_{L^{\infty}(\Omega)}
\le
C_1'\|\delta\mu_a\|_{L^{\infty}(\Omega)}\|u_0\|_{L^{\infty}(\Omega)}
+
C_1'\|\delta\mu_a\|_{L^{\infty}(\Omega)}\|u-u_R\|_{L^{\infty}(\Omega)},
\ee
where $C_1'$ is a positive constant. By moving the last term on the right-hand side to the left-hand side, we arrive at the following inequality for sufficiently small $\|\delta\mu_a\|_{L^{\infty}(\Omega)}$.
\be
\|u-u_R\|_{L^{\infty}(\Omega)}\le C_1\|\delta\mu_a\|_{L^{\infty}(\Omega)}\|u_0\|_{L^{\infty}(\Omega)},
\label{lem1}
\ee
where $C_1$ is a positive constant.

The outgoing light is detected at $\vv{r}_d^{(p)}\in\pp\Omega$. Let us introduce the data $\phi^{(p)}$ as
\be
\phi^{(p)}=\ln\frac{u_0(\vv{r}_d^{(p)})}{u(\vv{r}_d^{(p)})}.
\ee
In the first Rytov approximation we have $\phi^{(p)}\approx\phi_{R2}^{(p)}\approx\phi_R^{(p)}$, where
\ba
\phi_R^{(p)}
&=
-\frac{v_1(\vv{r}_d^{(p)})}{u_0(\vv{r}_d^{(p)})},
\\
\phi_{R2}^{(p)}
&=
-\frac{v_1(\vv{r}_d^{(p)})}{u_0(\vv{r}_d^{(p)})}
+\frac{1}{2}\left(\frac{v_1(\vv{r}_d^{(p)})}{u_0(\vv{r}_d^{(p)})}\right)^2
-\frac{v_2(\vv{r}_d^{(p)})}{u_0(\vv{r}_d^{(p)})}.
\ea

Note that $u_0,u_R$ are positive. Using (\ref{lem1}), we have
\ba
\left|e^{-\phi^{(p)}}\right|
&=
\left|\frac{u-u_R}{u_0}+\frac{u_R}{u_0}\right|
\\
&\le
\frac{C_1\|\delta\mu_a\|_{L^{\infty}(\Omega)}\|u_0\|_{L^{\infty}(\Omega)}}{u_0}+e^{v_1/u_0}.
\ea
The above inequality implies that $\phi^{(p)}\to0$ as $\|\delta\mu_a\|_{L^{\infty}(\Omega)}\to 0$.

\section{Discretization}
\label{discr}

Let $M$ be an even integer. We introduce a {\em spin} variable of discrete values as
\be
S(\vv{r})=0,\pm1,\pm2,\dots,\pm\frac{M}{2},\quad\vv{r}\in\Omega_{\rm ROI}.
\ee
To explain how $S(\vv{r})$ can be related to $\delta\mu_a(\vv{r})$, we will consider the following two cases. In both cases, the two-dimensional space ($d=2$) is assumed.

\subsection{Single-spin model}

Let us assume
\be
\delta\mu_a(\vv{r})=\eta f_{\bar{a}}(x)\delta(y-y_0),
\label{fwd:dela}
\ee
where $\eta,y_0$ are given positive constants. In this case, $\Omega_{\rm ROI}$ is a point $(x,y)=(0,y_0)$. Here, we assume that $f_{\bar{a}}(x)$ is given by
\be
f_{\bar{a}}(x)=
\left[\bar{a}^3+3\left(1+\frac{\tanh{x^2}}{10}\right)\bar{a}^2\right]
\left(1-\tanh{x^2}\right),
\label{fwd:delaf}
\ee
where $\bar{a}$ is a constant. Thus, $\delta\mu_a$ is determined by $\bar{a}$. This $\bar{a}$ is the unknown parameter to be reconstructed. In the forward data.

When $\bar{a}$ is sought, we assume that the minimum $a^{\rm(min)}$ and maximum $a^{\rm(max)}$ of candidates of $\bar{a}$ are known a priori. Let $a$ be a candidate of $\bar{a}$. We give $a$ as
\be
a=a^{\rm(min)}+\frac{a^{\rm(max)}-a^{\rm(min)}}{M+1}\left(S+\frac{M}{2}+1\right)
\in(a^{\rm(min)},a^{\rm(max)}].
\ee

\subsection{Multi-spin model}

Next we suppose that the support of $\delta\mu_a$ is unknown but we know $\delta\mu_a$ is zero outside a domain $\Omega_{\rm ROI}$ ($\mathop{\mathrm{supp}}\delta\mu_a\subset\Omega_{\rm ROI}$).

We divide the region of interest $\Omega_{\rm ROI}$ into cells $\omega_i\subset\Omega$ ($i=1,\dots,N$) such that $\Omega_{\rm ROI}=\bigcup_{i=1}^N\omega_i$ and $\omega_{i_1}\cap\omega_{i_2}=\varnothing$ if $i_1\neq i_2$. Let $|\omega|$ denote the area (volume) of subdomains $\omega_i$. We suppose $N>M_{\rm SD}$, which means the inverse problem is underdetermined. Let us assume $\delta\mu_a\in[0,\delta\mu_a^{\rm(max)}]$ with a constant $\delta\mu_a^{\rm(max)}>0$. We introduce
\be
\delta\mu_a(\vv{r}_i)=\delta\mu_a^{\rm(max)}\left(\frac{S_i}{M}+\frac{1}{2}\right),\quad
S_i=0,\pm1,\pm2,\dots,\pm\frac{M}{2}
\label{spins}
\ee
for $i=1,\dots,N$. Here, $\vv{r}_i\in\omega_i$ is a representative point in $\omega_i$. Thus $\delta\mu_a$ is discretized in $[0,\delta\mu_a^{\rm(max)}]$ by $M+1$ values at each point $\vv{r}_i$.

\section{Spin Hamiltonian: Case 1}
\label{spinh1}

We consider the single-spin model (\ref{fwd:dela}) for the time-dependent diffusion equation (\ref{de2d}).

We note that the Green's function for (\ref{de2d}) is given by
\be
\begin{aligned}
G(\vv{r},t;\vv{r}',s)
&=
\frac{e^{-\bar{\mu}_ac(t-s)}}{4\pi D_0(t-s)}e^{-\frac{(x-x')^2}{4D_0c(t-s)}}
\Biggl[e^{-\frac{(y-y')^2}{4D_0c(t-s)}}+e^{-\frac{(y+y')^2}{4D_0c(t-s)}}
\\
&-
\frac{\sqrt{4\pi D_0c(t-s)}}{\ell}e^{-\frac{(y+y')^2}{4D_0c(t-s)}}
e^{\left(\frac{y+y'+2D_0c(t-s)/\ell}{\sqrt{4D_0c(t-s)}}\right)^2}
\mathrm{erfc}\left(\frac{y+y'+2D_0c(t-s)/\ell}{\sqrt{4D_0c(t-s)}}\right)
\Biggr]
\end{aligned}
\ee
for $t>s$, and otherwise $G(\vv{r},t;\vv{r}',s)=0$. Here, $\ell=D_0\zeta$ and the complementary error function is given by $\mathop{\mathrm{erfc}}(x)=(2/\sqrt{\pi})\int_x^{\infty}\exp(-t^2)\,dt$. Using the Green's function, we obtain
\be
u_0(\vv{r},t)=
\frac{g_0e^{-\bar{\mu}_act}}{2\pi D_0t}e^{-\frac{(x-x_s^i)^2+y^2}{4D_0ct}}
\left[1-
\frac{\sqrt{\pi D_0ct}}{\ell}e^{\left(\frac{y+2D_0ct/\ell}{\sqrt{4D_0ct}}\right)^2}
\mathop{\mathrm{erfc}}\left(\frac{y+2D_0ct/\ell}{\sqrt{4D_0ct}}\right)\right].
\ee

Let us define
\be
\begin{aligned}
g(y,t;y',s)&=
\frac{1}{4\pi D_0(t-s)}e^{-\frac{(y+y')^2}{4D_0c(t-s)}}
\Biggl[1+e^{\frac{(y+y')^2-(y-y')^2}{4D_0c(t-s)}}-
\frac{\sqrt{4\pi D_0c(t-s)}}{\ell}
\\
&\times
e^{\left(\frac{y+y'}{2\sqrt{D_0c(t-s)}}+\frac{\sqrt{D_0c(t-s)}}{\ell}\right)^2}
\mathrm{erfc}\left(\frac{y+y'}{2\sqrt{D_0c(t-s)}}+\frac{\sqrt{D_0c(t-s)}}{\ell}
\right)
\Biggr].
\end{aligned}
\ee
Then we have
\be
\int_0^tG(\vv{r},t;\vv{r}',s)u_0(\vv{r}',s)\,ds=
e^{-\bar{\mu}_act}\int_0^te^{-\frac{(x-x')^2}{4D_0c(t-s)}}
e^{-\frac{(x'-x_s^{(p)})^2}{4D_0cs}}g(y,t;y',s)g(y',s;0,0)\,ds.
\ee
Therefore we obtain
\be
\begin{aligned}
u_R(\vv{r},t;\vv{r}_s^{(p)})
&=
u_0(\vv{r},t;\vv{r}_s^{(p)})
\exp\Bigg[-
\frac{e^{-\bar{\mu}_act}}{u_0(\vv{r},t;\vv{r}_s^i)}
\int_0^{\infty}\int_{-\infty}^{\infty}\delta\mu_a(\vv{r}')
\\
&\times
\left(\int_0^te^{-\frac{(x-x')^2}{4D_0c(t-s)}}
e^{-\frac{(x'-x_s^{(p)})^2}{4D_0cs}}g(y,t;y',s)g(y',s;0,0)\,ds\right)\,dx'dy'
\Bigg].
\end{aligned}
\label{fwd:ufinal}
\ee
We will use (\ref{fwd:ufinal}) to obtain the forward data in our numerical experiment.

By substituting the form (\ref{fwd:dela}) for $\delta\mu_a$ in (\ref{fwd:ufinal}), we obtain
\be
\begin{aligned}
u_R(\vv{r}_d^{(p)},t;\vv{r}_s^{(p)};a)
&=
u(\vv{r}_d^{(p)},t)
\\
&=
u_0(\vv{r}_d^{(p)},t;\vv{r}_s^{(p)})
\exp\Bigg[-\frac{\eta e^{-\bar{\mu}_act}}{u_0(\vv{r}_d^{(p)},t;\vv{r}_s^{(p)})}\int_0^tg(0,t;y_0,s)
\\
&\times
g(y_0,s;0,0)\left(\int_{-\infty}^{\infty}f_a(x')
e^{-\frac{(x_d^{(p)}-x')^2}{4D_0c(t-s)}}e^{-\frac{(x'-x_s^{(p)})^2}{4D_0cs}}
\,dx'\right)\,ds\Bigg],
\end{aligned}
\label{fwd:formula}
\ee
where
\be
u_0=\frac{g_0e^{-\bar{\mu}_act}}{2\pi D_0t}e^{-\frac{(x_d^{(p)}-x_s^{(p)})^2}{4D_0ct}}
\left[1-\frac{\sqrt{\pi D_0ct}}{\ell}
e^{\left(\frac{\sqrt{D_0ct}}{\ell}\right)^2}
\mathrm{erfc}\left(\frac{\sqrt{D_0ct}}{\ell}\right)\right].
\ee
We obtain
\ba
\phi_R^{(p)}&=
\frac{\eta e^{-\bar{\mu}_act}}{u_0(\vv{r}_d^{(p)},t;\vv{r}_s^{(p)})}
\int_0^tg(0,t;y_0,s)g(y_0,s;0,0)
\\
&\times
\left(\int_{-\infty}^{\infty}f_a(x')
e^{-\frac{(x_d^{(p)}-x')^2}{4D_0c(t-s)}}e^{-\frac{(x'-x_s^{(p)})^2}{4D_0cs}}
\,dx'\right)\,ds.
\ea
Let us discretize time as $t\to t_j$ ($j=1,\dots,M_t$). With this discretization, we write $\tilde{\phi}_R^{(p)}$ instead of $\phi_R^{(p)}$.

In order to use SA instead of the naive Metropolis-Hastings Markov-chain Monte Carlo \cite{Jiang-etal19}, we treat the cost function as a spin Hamiltonian. We let $\Phi^{(p)}$ denote the experimentally obtained data corresponding to $\phi^{(p)}$. We can solve the inverse problem for our DOT by minimizing the cost function which is given by
\be
\mathcal{H}(S)=\frac{1}{2}\sum_{p=1}^{M_{\rm SD}}\sum_{j=1}^{M_t}\left|\Phi^{(p)}-\tilde{\phi}_R^{(p)}\right|^2+\alpha|S-S^{(0)}|,
\label{hami1}
\ee
where $\alpha$ is the regularization parameter and $S^{(0)}$ is an initial guess. We put $\alpha=0$.

The cost function in (\ref{hami1}) has one local minimum and one global minimum \cite{Jiang-etal19}. To see this structure of $\mathcal{H}$, we introduce
\be
h(t)=\frac{1}{t}\exp\left(-\frac{y_0^2}{4Dct}\right)
\left[1-\frac{\sqrt{\pi Dct}}{\ell}
e^{\left(\frac{y_0}{2\sqrt{Dct}}+\frac{\sqrt{Dct}}{\ell}\right)^2}
\mathrm{erfc}\left(\frac{y_0}{2\sqrt{Dct}}+\frac{\sqrt{Dct}}{\ell}\right)
\right].
\ee
The following form is obtained using Eq.~(\ref{fwd:formula}). By neglecting noise, we have
\be
\begin{aligned}
&
u\left(\vv{r}_d^{(p)},t_j;\vv{r}_s^{(p)};\bar{a}\right)-
u\left(\vv{r}_d^{(p)},t_j;\vv{r}_s^{(p)};a\right)
\\
&=
\frac{\eta e^{-\mu_{a0}ct_j}}{(2\pi D)^2}\int_0^{t_j}h(t_j-s)h(s)
\\
&\times
\left[\int_{-\infty}^{\infty}d_{\bar{a}}(a,x')\left(1-\tanh{{x'}^2}\right)
e^{-\frac{(x_d^{(p)}-x')^2}{4Dc(t_j-s)}}e^{-\frac{(x'-x_s^{(p)})^2}{4Dcs}}
\,dx'\right]\,ds,
\end{aligned}
\ee
where $d_{\bar{a}}(a,x')=\xi(\bar{a},x')-\xi(a,x')$ with
\be
\xi(a,x')=a^2\left[a+3\left(1+\frac{\tanh{{x'}^2}}{10}\right)\right].
\ee
For a given $x'$, the function $|d_{\bar{a}}(a;x')|^2$ has the global minimum at $a=\bar{a}>0$, a local minimum at $a=-2\left(1+\frac{\tanh{{x'}^2}}{10}\right)$<0, and a local maximum at $a=0$. This structure of $|d_{\bar{a}}(a;x')|^2$ implies that iterative methods, which always look for a position that lowers the value of the cost function, do not work when the initial guess $a_0$ is negative. 

\section{Spin Hamiltonian: Case 2}
\label{spinh}

Here we consider the multi-spin model (\ref{spins}) for the time-independent diffusion equation (\ref{de0}).

Let us define
\be
K_{p,i}=
\delta\mu_a^{\rm(max)}|\omega|\frac{G(\vv{r}_d^{(p)},\vv{r}_i)G(\vv{r}_i,\vv{r}_s^{(p)})}{G(\vv{r}_d^{(p)},\vv{r}_s^{(p)})},
\ee
where $G$ is the Green's function for (\ref{de0}) with $\mu_a=\bar{\mu}_a$. Let $\tilde{\phi}_R^{(p)}$, $\tilde{\phi}_{R2}^{(p)}$ be discretized versions of $\phi_R^{(p)}$, $\phi_{R2}^{(p)}$, respectively. Assuming that $\omega_i$ is small, we have
\be
\phi_R^{(p)}\approx\tilde{\phi}_R^{(p)},\quad
\phi_{R2}^{(p)}\approx\tilde{\phi}_{R2}^{(p)},
\ee
where
\be
\tilde{\phi}_R^{(p)}=
\sum_{i=1}^NK_{p,i}\left(\frac{S_i}{M}+\frac{1}{2}\right)
\ee
and
\ba
\tilde{\phi}_{R2}^{(p)}
&=
\sum_{i=1}^NK_{p,i}\left(\frac{S_i}{M}+\frac{1}{2}\right)
\\
&+
\frac{1}{2}\sum_{i_1=1}^N\sum_{i_2=1}^NK_{p,i_1}K_{p,i_2}
\left(\frac{S_{i_1}}{M}+\frac{1}{2}\right)
\left(\frac{S_{i_2}}{M}+\frac{1}{2}\right)
\\
&-
|\omega|^2\sum_{i_1=1}^N\sum_{i_2=1}^N\left[G(\vv{r}_d^{(p)},\vv{r}_s^{(p)})\right]^{-1}
\\
&\times
G(\vv{r}_d^{(p)},\vv{r}_{i_1})\delta\mu_a(\vv{r}_{i_1})G(\vv{r}_{i_1},\vv{r}_{i_2})\delta\mu_a(\vv{r}_{i_2})G(\vv{r}_{i_2},\vv{r}_s^{(p)}).
\ea

Let $\vv{S}=(S_i)$ ($i=1,\dots,N$) denote the spin configuration. We will look for $\vv{S}^*$ which minimizes the following cost function $\Psi(\vv{S})$.
\be
\Psi(\vv{S})=\frac{1}{2}\sum_{p=1}^{M_{\rm SD}}\left|\Phi^{(p)}-\tilde{\phi}_{R2}^{(p)}\right|^2+\alpha\sum_{i=1}^N|S_i-S_i^{(0)}|,
\ee
where $\alpha>0$ is the regularization parameter and $\vv{S}^{(0)}=(S_i^{(0)})\in\Rm^N$ is an initial guess, which we set $S_i^{(0)}=-M/2$ ($i=1,\dots,N$).

To compare terms of order $(\delta\mu_a^{\rm(max)})^2$, let us introduce
\ba
g_1&=K_{p,i_1}K_{p,i_2},
\\
g_2&=g_1-
2\left(\delta\mu_a^{\rm(max)}\right)^2|\omega|^2\frac{G(\vv{r}_d^{(p)},\vv{r}_{i_1})G(\vv{r}_{i_1},\vv{r}_{i_2})G(\vv{r}_{i_2},\vv{r}_s^{(p)})}{G(\vv{r}_d^{(p)},\vv{r}_s^{(p)})}.
\ea
We have
\ba
\left|\phi^{(p)}g_2\right|
&=
|g_1|\left|\phi^{(p)}\frac{g_2}{g_1}\right|
\\
&=
|g_1|\left|\phi^{(p)}\right|\left|
1-2\frac{G(\vv{r}_d^{(p)},\vv{r}_s^{(p)})G(\vv{r}_{i_1},\vv{r}_{i_2})}{G(\vv{r}_d^{(p)},\vv{r}_{i_2})G(\vv{r}_{i_1},\vv{r}_s^{(p)})}\right|.
\ea
Recall that $|\phi^{(p)}|\to0$ as $\|\delta\mu_a\|_{L^{\infty}(\Omega)}\to0$. Thus for arbitrary $p,i_1,i_2$ ($p=1,\dots,M_{\rm SD}$, $i_1,i_2=1,\dots,N$),
\be
|g_1|\ge|\phi^{(p)}g_2|
\label{thm1}
\ee
when $\|\delta\mu_a\|_{L^{\infty}(\Omega)}$ is sufficiently small.

By neglecting $g_2$, we have $\Psi=\mathcal{H}(\vv{S})+\mbox{const.}$, where only the first term $\mathcal{H}(\vv{S})$ on the right-hand side depends on $S_i$. The Hamiltonian $\mathcal{H}(\vv{S})$ is given by
\be
\mathcal{H}(\vv{S})=-\sum_{i=1}^N\sum_{j=1}^NJ_{ij}S_iS_j-\sum_{i=1}^Nh_iS_i,
\label{Ham}
\ee
where
\ba
J_{ij}
&=
\frac{-1}{2M^2}\sum_{p=1}^{M_{\rm SD}}K_{p,i}K_{p,j},
\\
h_i
&=
M\sum_{j=1}^NJ_{ij}+
\frac{1}{M}\left(\sum_{p=1}^{M_{\rm SD}}\Phi^{(p)}K_{p,i}-\alpha\right).
\ea
The spin interactions are symmetric: $J_{ij}=J_{ji}$. We note that the modulus $|h_i|$ of the magnetic field becomes large when the hyperparameter $\alpha$ in the regularization term is large. In this case, the spin Hamiltonian has a unique ground state ($\vv{S}^*=\vv{S}^{(0)}$) for sufficiently large $\alpha$.

Thus our optical tomography is reformulated as the problem of searching the ground state of $\mathcal{H}(\vv{S})$, i.e., the spin configuration $\vv{S}^*=\argmin_{S}\mathcal{H}(\vv{S})$.

\section{Simulated annealing}
\label{SA}

The algorithm of the simulated annealing is described below. Although we assume Case 2 of multiple spins, the single spin in Case 1 can be similarly implemented.

The Hamiltonian $\mathcal{H}(\vv{S})$ in (\ref{Ham}) can be seen as a Hamiltonian of spins that interact with each other in a solid state material. In a real physical spin system, spins fluctuate depending on temperature and the spin system reaches its lowest energy state if temperature gradually decreases. In SA, we seek the optimal solution by mimicking this physical process with thermal fluctuations. To this end, we introduce the parameter $T$ which corresponds to temperature. We will find the ground state $\vv{S}^*$ of $\mathcal{H}(\vv{S})$ in (\ref{Ham}) by gradually decreasing $T$ from $T_{\rm high}$ to $T_{\rm low}$.

The partition function $Z$ is given by
\be
Z=\sum_{\{S_i\}}e^{-\beta\mathcal{H}(\vv{S})},
\ee
where $\beta=1/T$ is the inverse temperature. Here, we used the notation $\sum_{\{S_i\}}=\sum_{S_1}\cdots\sum_{S_N}$. We treat $\{S_i\}=\{S_1,S_2,\dots,S_N\}$ as random variables. We give the probability density function $\pi(\vv{S})$ as the Boltzmann distribution given by
\be
\pi(\vv{S})=\frac{e^{-\beta\mathcal{H}(\vv{S})}}{Z}.
\label{probS}
\ee
We see that $\pi(\vv{S})$ is large if $\mathcal{H}(\vv{S})$ is small. When the temperature is sufficiently low, i.e., $\beta$ is large, $\pi(\vv{S}^*)$ becomes significantly larger than the probability of other configurations.

Although calculating the denominator $Z$ on the right-hand side of (\ref{probS}) is difficult, we can find $\vv{S}^*$ by using the Metropolis algorithm \cite{Metropolis-etal53,Landau-Binder}. The proposal distribution $q(S_i'|\vv{S})$ is given for each spin, say the $i$th spin, and the value of $S_i'$ is generated with equal probability. We note that
\ba
&
\frac{\pi(\vv{S}')}{\pi(\vv{S})}=e^{\beta(\mathcal{H}(\vv{S})-\mathcal{H}(\vv{S}'))}
\\
&=
\exp\left[\beta\left(J_{ii}\left({S_i'}^2-S_i^2\right)+
\left(2\sum_{j=1\atop j\neq i}^NJ_{ij}S_j+h_i\right)\left(S_i'-S_i\right)\right)\right].
\ea

While the naive Metropolis-Hastings algorithm does not work in high dimensions \cite{Au-Beck01,Katafygiotis-Zuev08}, the Metropolis algorithm or single-component Metropolis-Hastings algorithm can handle many unknown parameters \cite{Gilks-Richardson-Spiegelhalter}. Although we use the standard Metropolis algorithm as the first attempt of the single-component Metropolis-Hastings algorithm for DOT, different algorithms with faster relaxation have been developed for spin systems \cite{Landau-Binder}. Our algorithm is summarized as the following six steps.

\vspace{1em}
\begin{enumerate}
\setlength\itemsep{1em}
\item
Start with a small $\beta=1/T_{\rm high}>0$. Give $(S_i)$ randomly as an initial guess. Then set $i=1$.
\item
Compute $h_{{\rm eff},i}=2\sum_{j=1\;(j\neq i)}^NJ_{ij}S_j+h_i$.
\item
Calculate $w=-\beta[h_{{\rm eff},i}(S_i'-S_i)+J_{ii}({S_i'}^2-S_i^2)]$, where $S_i'\sim q(\cdot|\vv{S})$ is randomly chosen.
\item
Set $S_i=S_i'$ if $w\le0$. Otherwise put $S_i=S_i'$ with probability $e^{-w}$.
\item
Set $i=1$ if $i=N$. Otherwise set $i=i+1$. Return to Step 2. After repeating several loops from Step 2 to Step 5 until the initial large fluctuation ceases, proceed to Step 6.
\item
Decrease temperature and go to Step 2. If the temperature reaches $T_{\rm low}$, finish the iteration.
\end{enumerate}

In this paper, we decrease $T$ as
\be
T-10^{{\rm int}(\log_{10}T)-2}\quad\rightarrow\quad T.
\label{Tplan}
\ee

\section{Numerical test}
\label{numerics}

\subsection{Single spin}

We take the unit of length and unit of time to be ${\rm mm}$ and ${\rm ps}$, respectively. On the $x$-axis, we place two sources at $(x,y)=(-20,0)$, $(20,0)$ and three detectors at $(x,y)=(-40,0)$, $(0,0)$, $(40,0)$. As a result we have $M_{\rm SD}=6$. We set $D_0=0.33$, $\mu_a=0.02$, $\mathfrak{n}=1.37$. Suppose that there is absorption inhomogeneity at depth $5$. For $\delta\mu_a$, we put $\eta=300/c$, $y_0=5$, and
\be
\bar{a}=1.5.
\ee
When we prepare $\Phi^{(p)}$, we added $3\%$ Gaussian noise. We set
\be
t_j=j\Delta_t\quad(j=1,\dots,M_t),\quad \Delta_t=5,\quad M_t=500.
\ee
Furthermore we set $M=512$, $a^{\rm(min)}=-3$, and $a^{\rm(max)}=3$.

\begin{figure}[ht!]
\centering
\includegraphics[width=0.8\textwidth]{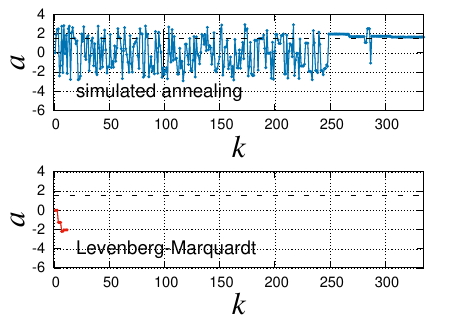}
\caption{
(Top) The value $a$ is plotted against Monte Carlo steps for the simulated annealing in Sec.~\ref{SA}. (Bottom) The reconstruction is done with the Levenberg–Marquardt algorithm.
}
\label{toy:fig1}
\end{figure}

It is known that iterative methods cannot reach $\bar{a}=1.5$ when the initial value is negative \cite{Jiang-etal19}. We set the initial value $a_0$ of $a$ as
\be
a_0=-0.01.
\ee
In Fig.~\ref{toy:fig1}, we compare the simulated annealing developed in this paper with the Levenberg–Marquardt algorithm \cite{Fletcher71,Levenberg44,Marquardt63}, which is one of iterative methods. The Levenberg-Marquardt algorithm fails to arrive at the correct answer, whereas SA converges to the correct value. In our simulation, converged values for the Levenberg-Marquardt algorithm and SA are $-2.05$ and $1.68$, respectively. As is mentioned in the end of Sec.~\ref{spinh1}, the cost function has a local minimum at a negative value of $a$. The calculation of the Levenberg-Marquardt algorithm is trapped by this local minimum. On the other hand, $a$ goes back and forth between the local minimum and global minimum in SA, and then falls in the global minimum as temperature decreases.

\subsection{Multiple spins}

Let us  consider DOT in the half space $\Omega$ , i.e., $\Omega=\{\vv{r}\in\Rm^2;\;-\infty<x<\infty,\;0<y<\infty\}$. The boundary $\pp\Omega$ is the $x$-axis. Measurements are performed on the $x$-axis: $\vv{r}_s^{(p)}=(x_s^{(p)},0^+)$, $\vv{r}_d^{(p)}=(x_d^{(p)},0)$. In this case, the Green's function is explicitly given by
\be
G(\vv{r},\vv{r}')=
\frac{1}{2\pi D_0}\int_0^{\infty}\frac{\cos(q(x-x'))}{\lambda}
\left(
e^{-\lambda|y-y'|}+\frac{\ell\lambda-1}{\ell\lambda+1}e^{-\lambda(y+y')}
\right)\,dq,
\label{green}
\ee
where $\ell=\zeta D_0$ and $\lambda=\lambda(q)=\sqrt{\frac{\bar{\mu}_a}{D_0}+q^2}$. Equation (\ref{green}) implies $\|G(\vv{r},\cdot)\|_{L^1(\Omega)}<\infty$ and $\|G(\vv{r},\cdot)\|_{L^{\infty}(\Omega)}<\infty$ for any $\vv{r}\in\Omega$. The integrand in (\ref{green}) decays slowly when $y,y'$ are small and is oscillatory. The numerical integration in (\ref{green}) can be done by the double-exponential formula \cite{Ooura-Mori91}.

We use $M_{\rm SD}=240$ source-detector pairs from $16$ sources and $15$ detectors:
\ba
x_s^{(p)}
&=\pm2,\pm6,\dots,\pm30\,{\rm mm},
\nonumber \\
x_d^{(p)}
&=0,\pm4,\pm8,\dots,\pm28\,{\rm mm}.
\ea
In addition to the refractive index $\mathfrak{n}=1.37$, we set
\be
\bar{\mu}_a=0.02\,{\rm mm}^{-1},\quad D_0=0.33\,{\rm mm}^{-1}.
\ee
We assume absorption inhomogeneity of the shape of disks in the medium. Inside the disks, we set
\be
\delta\mu_a=0.2\,{\rm mm}^{-1}.
\ee
We compute the forward data $\Phi^{(p)}$ ($p=1,\dots,M_{\rm SD}$) by the finite-difference scheme with the Gauss-Seidel method. We added $3\%$ Gaussian noise to $u(\vv{r}_d^{(p)})$ and $u_0(\vv{r}_d^{(p)})$. Assuming the diffuse surface reflection, we have \cite{Egan-Hilgeman}
\be
\zeta=2\frac{1+r_d}{1-r_d},
\ee
where
\be
r_d=-1.4399\mathfrak{n}^{-2}+0.7099\mathfrak{n}^{-1}+0.6681+0.0636\mathfrak{n}.
\ee

We take $(2N_x+1)N_y=1830$ cells ($N_x=30$, $N_y=30$) of the area $|\omega|=h^2$ ($h=1\,{\rm mm}$). Moreover, the following parameter values were used for SA.
\be
\alpha=0.01,\quad T_{\rm high}=10^{-5},\quad T_{\rm low}=10^{-10},\quad M=256.
\ee

\begin{figure}[ht!]
\centering
\includegraphics[width=0.4\textwidth]{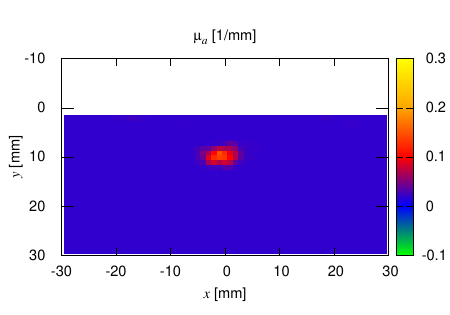}
\hspace{5mm}
\includegraphics[width=0.4\textwidth]{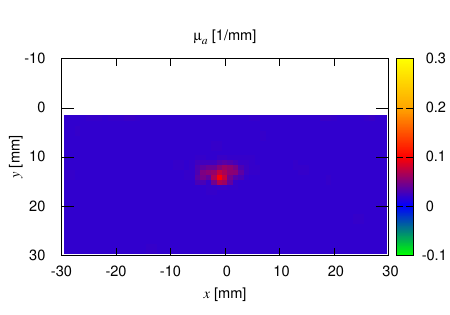}
\caption{
(Left) One absorber at the depth $10\,{\rm mm}$. (Right) One absorber at a deeper position at $y=15\,{\rm mm}$.
}
\label{sa:fig1}
\end{figure}

\begin{figure}[ht!]
\centering
\includegraphics[width=0.4\textwidth]{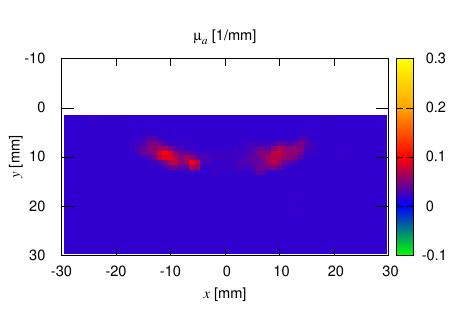}
\caption{
Two absorbers with separation $20\,{\rm mm}$ reconstructed by the proposed method.
}
\label{sa:fig2}
\end{figure}

Figures \ref{sa:fig1} and \ref{sa:fig2} show reconstructed images by the proposed method. In Fig.~\ref{sa:fig1}, a disk target of radius $2.5\,{\rm mm}$ was placed at the depth (Left) $10\,{\rm mm}$ and (Right) $15\,{\rm mm}$, that is, the $x$ and $y$ coordinates of the disk center are $x=0\,{\rm mm}$, and $y=10\,{\rm mm}$ or $15\,{\rm mm}$. In general, the reconstruction at deeper positions is more difficult because the signal becomes noisy and the relative contribution of the regularization term in the cost function becomes larger \cite{Shimokawa-etal12}. As is expected, the reconstruction at a deeper position is more difficult. 

\begin{figure}[ht!]
\centering
\includegraphics[width=0.4\textwidth]{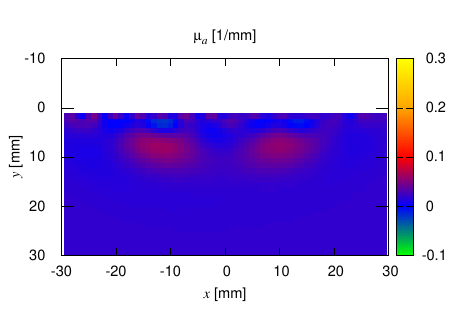}
\hspace{5mm}
\includegraphics[width=0.4\textwidth]{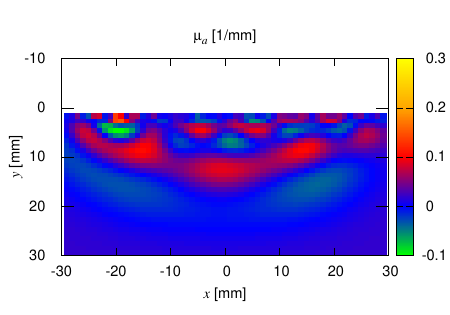}
\caption{
Two absorbers with separation $20\,{\rm mm}$. (Left) The reconstruction by the truncated SVD with $52$ singular values. (Right) The reconstruction by the truncated SVD with $80$ singular values.
}
\label{sa:fig3}
\end{figure}

Next we consider two disks of radius $2.5\,{\rm mm}$ in Fig.~\ref{sa:fig2}. The centers of the disks are placed at $(x,y)=(\pm10\,{\rm mm},\,10\,{\rm mm})$. Figure \ref{sa:fig2} shows the reconstruction by the proposed scheme using the same parameters described above for Fig.~\ref{sa:fig1}. For comparison, we also used the truncated SVD. Reconstructions from the truncated SVD are presented in Fig.~\ref{sa:fig3} with (Left) $52$ largest singular values and (Right) $80$ largest singular values, respectively. In the reconstructed images by the truncated SVD, the separation between the two targets is less sharp.

In each numerical calculation, there are $5\times 10^5$ unknown parameter values. The computation time to find $\vv{S}^*$ was about three minutes on a laptop computer (MacBook Pro with 2.3 GHz Intel Core i5 and 8 GB memory).

\section{Concluding remarks}
\label{concl}

We have developed a stochastic approach of DOT using the Metropolis algorithm of the Monte Carlo method. Due to the fact that the inverse problem for DOT is severely ill-posed, the reconstructed values are always blurred at least to some extent. This motivated us to introduce spins, which have discrete values. The effect of the discretization of the parameter is subtle for the single-spin model but becomes important when the number of spins increases in the multi-spin model because the discretization reduces the search space.

For linearized inverse problems, the truncated SVD is often used \cite{Arridge99,Arridge-Schotland09}. Although the truncated SVD has a filtering property similar to the Tikhonov regularization, the penalty term of the proposed method is not restricted to the $2$-norm. In this paper, the $1$-norm was used. This might explain the difference of the quality of reconstruction in Fig.~\ref{sa:fig3}. Another advantage of our method is that it is feasible to use a priori information and set the lower and upper bounds of $\delta\mu_a$, whereas the reconstructed $\delta\mu_a$ by the truncated SVD is unbounded. This feature is reflected in the difference between the proposed method and SVD in Fig.~\ref{sa:fig3}. As is seen for the single-spin model, the proposed method can treat both linear and nonlinear inverse problems.

The single-spin model provides a nonlinear inverse problem. Although it is a linearized inverse problem, the multi-spin model demonstrates that a large number of unknowns can be treated by our numerical scheme. It is a natural next step to consider nonlinear inverse problems with multiple spins. For this, nonlinear terms such as $u_{R2}$ in the Rytov series have to be taken into account. When the landscape has more complicated structure, it is important to decrease temperature gradually.

To achieve fast convergence, it is crucial to bring the spin system efficiently to the thermal equilibrium state. Various Monte Carlo methods for spin systems have been developed to reduce the computation time. These techniques may be implemented in our approach, which solves the inverse problem of DOT from the viewpoint of statistical mechanics of spin systems. Existing Monte Carlo methods for spin systems include the following algorithms. The replica exchange method prepares multiple copies of a spin system to let the system at low temperatures escape from local minima \cite{Hukushima-Nemoto96}. To minimize the average rejection rate, the MCMC algorithm by Suwa and Todo breaks the detailed balance condition while preserving the balance condition \cite{Suwa-Todo10}. Quantum annealing makes use of quantum fluctuations of spins whereas SA uses thermal fluctuations \cite{Kadowaki-Nishimori98}. Quantum annealing also has potential to provide an efficient algorithm for DOT.

\section*{Acknowledgments}
YJ is supported by the National Natural Science Foundation of China (No.~11971121). MM is supported by JSPS KAKENHI Grant Numbers JP17K05572, JP18K03438 and by HUSM Grant-in-Aid. NT is supported by JSPS KAKENHI Grant Number JP16K05418.


\bibliography{myref}

\begin{thebibliography}{10}

\bibitem{Arridge99}
S.~R. Arridge.
\newblock Optical tomography in medical imaging.
\newblock {\em Inverse Problems}, 15:R41--R93, 1999.

\bibitem{Arridge-etal06}
S.~R. Arridge, J.~P. Kaipio, V.~Kolehmainen, M.~Schweiger, E.~Somersalo,
  T.~Tarvainen, and M.~Vauhkonen.
\newblock Approximation errors and model reduction with an application in
  optical diffusion tomography.
\newblock {\em Inverse Problems}, 22:175--195, 2006.

\bibitem{Arridge-Schotland09}
S.~R. Arridge and J.~C. Schotland.
\newblock Optical tomography: forward and inverse problems.
\newblock {\em Inverse Problems}, 25:123010, 2009.

\bibitem{Au-Beck01}
S.-K. Au and J.~L. Beck.
\newblock Estimation of small failure probabilities in high dimensions by
  subset simulation.
\newblock {\em Probabilistic Engineering Mechanics}, 16:263--277, 2001.

\bibitem{Bal-etal13}
G.~Bal, I.~Langmore, and Y.~Marzouk.
\newblock Bayesian inverse problems with {M}onte {C}arlo forward models.
\newblock {\em Inv. Probl. Imag.}, 7:81--105, 2013.

\bibitem{Choe-etal05}
R.~Choe, A.~Corlu, K.~Lee, T.~Durduran, S.~D. Konecky, M.~Grosicka-Koptyra,
  S.~R. Arridge, B.~J. Czerniecki, D.~L. Fraker, A.~DeMichele, B.~Chance, M.~A.
  Rosen, and A.~G. Yodh.
\newblock Diffuse optical tomography of breast cancer during neoadjuvant
  chemotherapy: {A} case study with comparison to {MRI}.
\newblock {\em Med. Phys.}, 32:1128--1139, 2005.

\bibitem{Egan-Hilgeman}
W.~G. Egan and T.~W. Hilgeman.
\newblock {\em Optical Properties of Inhomogeneous Materials}.
\newblock Academic Press, 1979.

\bibitem{Fletcher71}
R.~Fletcher.
\newblock {\em A modified marquardt subroutine for nonlinear least squares
  (Report AERE-R 6799)}.
\newblock The Atomic Energy Research Establishment: Harwell, 1971.

\bibitem{Gibson-etal05}
A.~P. Gibson, J.~C. Hebden, and S.~R. Arridge.
\newblock Recent advances in diffuse optical imaging.
\newblock {\em Phys. Med. Biol.}, 50:R1--R43, 2005.

\bibitem{Gilks-Richardson-Spiegelhalter}
W.~R. Gilks, S.~Richardson, and D.~J. Spiegelhalter.
\newblock {\em Markov Chain Monte Carlo in Practice}.
\newblock Chapman and Hall, 1996.

\bibitem{Hukushima-Nemoto96}
K.~Hukushima and K.~Nemoto.
\newblock Exchange {M}onte {C}arlo method and application to spin glass
  simulations.
\newblock {\em J. Phys. Soc. Jpn.}, 65:1604--1608, 1996.

\bibitem{Jiang-etal19}
Y.~Jiang, Y.~Hoshi, M.~Machida, and G.~Nakamura.
\newblock A hybrid inversion scheme combining {M}arkov chain {M}onte {C}arlo
  and iterative methods for determining optical properties of random media.
\newblock {\em Appl. Sci.}, 9:3500, 2019.

\bibitem{Kadowaki-Nishimori98}
T.~Kadowaki and H.~Nishimori.
\newblock Quantum annealing in the transverse {I}sing model.
\newblock {\em Phys. Rev. E}, 58:5355--5363, 1998.

\bibitem{Katafygiotis-Zuev08}
L.~S. Katafygiotis and K.~M. Zuev.
\newblock Geometric insight into the challenges of solving high-dimensional
  reliability problems.
\newblock {\em Probabilistic Engineering Mechanics}, 23:208--218, 2008.

\bibitem{Keller69}
J.~B. Keller.
\newblock Accuracy and validity of the {B}orn and {R}ytov approximations.
\newblock {\em J. Opt. Soc. Am.}, 59:1003--1004, 1969.

\bibitem{Kirkinis08}
E.~Kirkinis.
\newblock Renormalization group interpretation of the {B}orn and {R}ytov
  approximations.
\newblock {\em J. Opt. Soc. Am. A}, 25:2499--2508, 2008.

\bibitem{Landau-Binder}
D.~P. Landau and K.~Binder.
\newblock {\em A Guide to Monte Carlo Simulations in Statistical Physics}.
\newblock Cambridge University Press, 2000.

\bibitem{Langmore13}
I.~Langmore, A.~B. Davis, and G.~Bal.
\newblock Multipixel retrieval of structural and optical parameters in a 2-{D}
  scene with a path-recycling {M}onte {C}arlo forward model and a new
  {B}ayesian inference engine.
\newblock {\em IEEE Trans. Geosci. Remote Sens.}, 51:2903--2919, 2013.

\bibitem{Levenberg44}
K.~Levenberg.
\newblock A method for the solution of certain non-linear problems in least
  squares.
\newblock {\em Quarterly Appl. Math.}, 2:164--168, 1944.

\bibitem{Markel-Schotland04}
V.~A. Markel and J.~C. Schotland.
\newblock Symmetries, inversion formulas, and image reconstruction for optical
  tomography.
\newblock {\em Phys. Rev. E}, 70:056616, 2004.

\bibitem{Markel-Schotland07}
V.~A. Markel and J.~C. Schotland.
\newblock On the convergence of the {B}orn series in optical tomography with
  diffuse light.
\newblock {\em Inverse Problems}, 23:1445--1465, 2007.

\bibitem{Marquardt63}
D.~W. Marquardt.
\newblock An algorithm for least-squares estimation of nonlinear parameters.
\newblock {\em SIAM J. Appl. Math.}, 11:431--441, 1963.

\bibitem{Metropolis-etal53}
N.~Metropolis, A.~W. Rosenbluth, M.~N. Rosenbluth, A.~H. Teller, and E.~Teller.
\newblock Equation of state calculations by fast computing machines.
\newblock {\em J. Chem. Phys.}, 21:1087--1092, 1953.

\bibitem{Ooura-Mori91}
T.~Ooura and M.~Mori.
\newblock The double exponential formula for oscillatory functions over the
  half infinite interval.
\newblock {\em J. Comput. Appl. Math.}, 38:353--360, 1991.

\bibitem{Schweiger-Arridge99}
M.~Schweiger and S.~R. Arridge.
\newblock Application of temporal filters to time resolved data in optical
  tomography.
\newblock {\em Phys. Med. Biol.}, 44:1699--1717, 1999.

\bibitem{Shimokawa-etal12}
T.~Shimokawa, T.~Kosaka, O.~Yamashita, N.~Hiroe, T.~Amita, Y.~Inoue, and
  M.~Sato.
\newblock Hierarchical {B}ayesian estimation improves depth accuracy and
  spatial resolution of diffuse optical tomography.
\newblock {\em Opt. Exp.}, 20:20427--20446, 2012.

\bibitem{Suwa-Todo10}
H.~Suwa and S.~Todo.
\newblock Markov chain {M}onte {C}arlo method without detailed balance.
\newblock {\em Phys. Rev. Lett.}, 105:120603, 2010.

\bibitem{Yamada-Okawa14}
Y.~Yamada and S.~Okawa.
\newblock Diffuse optical tomography: {P}resent status and its future.
\newblock {\em Opt. Rev.}, 21:185--205, 2014.

\bibitem{Zhao-Cooper17}
H.~Zhao and R.~J. Cooper.
\newblock Review of recent progress toward a fiberless, whole-scalp diffuse
  optical tomography system.
\newblock {\em Neurophotonics}, 5:011012, 2017.

\bibitem{Zhu-Poplack20}
Q.~Zhu and S.~Poplack.
\newblock A review of optical breast imaging: {M}ulti-modality systems for
  breast cancer diagnosis.
\newblock {\em European J. Radiology}, 129:109067, 2020.

\end{thebibliography}
\bibliographystyle{plain}

\end{document}